\documentstyle{marcdek}

\def\Bbb{\bf}

\def\R{{\Bbb R}}

\def\listit#1,#2{#1_1$, $#1_2,\ldots,$\ $#1_{#2}}

\def\normo#1{{\left\| #1 \right\|}}
\def\snormo#1{{\mathopen\| #1 \mathclose\|}}

\def\widedot{\,\cdot\,}
\def\normdot{\normo{\widedot}}

\def\modo#1{{\left| #1 \right|}}
\def\smodo#1{{\mathopen| #1 \mathclose|}}

\def\and{\mathop{\hbox{\ and\ }}}
\def\meas{\mathop{\hbox{measure}}}

\newcounter{rom}
\newenvironment{itemlist}{\begin{list}%
{\roman{rom})}{\usecounter{rom}%
\setlength\parsep{0in}\setlength\itemsep{0in}\setlength\topsep{0in}%
\setlength\partopsep{0in}}}%
{\end{list}}

\let\itemi\item
\let\itemii\item
\let\itemiii\item
\let\itemiv\item

\begin{document}

\title{The Hardy Operator and Boyd Indices}

\author{\uppercase{Stephen~J.~Montgomery-Smith}%
\thanks{Research supported in part by N.S.F.\ Grant
DMS 9201357.}\ \ \
Department of Mathematics, University of Missouri,
Columbia, Missouri 65211.}

\maketitle

\thispagestyle{empty}
\pagestyle{empty}

\section*{ABSTRACT}

We give necessary and sufficient conditions for the
Hardy operator to be bounded on a rearrangement invariant quasi-Banach
space in terms of its Boyd indices.

\section*{MAIN RESULTS}

A {\it rearrangement invariant space\/} $X$\ on $\R$\ is a 
set of measurable
functions (modulo functions equal almost everywhere) with a 
complete quasi-norm $\normdot_X$\
such that the following holds:
\begin{itemlist}
\itemi if $g^* \le f^*$\ and $f \in X$, then $g \in X$\ with
$\snormo g_X \le \snormo f_X$;
\itemii if $f$\ is simple with finite support then $f \in X$;
\itemiii either $f_n \searrow 0$\ implies $\snormo{f_n} \searrow 0$
\item[iii$'$)] or $0 \le f_n \nearrow f$\ and $\sup_n \snormo{f_n}_X < \infty$\
imply $f \in X$\ with $\snormo f_X = \sup_n \snormo{f_n}_X$.
\end{itemlist}

\noindent
Here $f^*$\ denotes the decreasing rearrangement of $\smodo f$, that is,
$f(s) = \sup\{ \,t : \meas\{ \smodo f > t\} > s \, \}$.

The properties of rearrangement invariant spaces that we will use will be
the {\it Boyd indices\/} defined as follows.  Given a number $0<a<\infty$,
we define the operator $D_a f(t) = f(at)$.  Then the lower Boyd index
of $X$\ is defined by
\begin{eqnarray*}
   p_X &=& \sup\{\, p : \exists c \ \forall a<1 \ \snormo{D_a f}_X \le c \,
   a^{-1/p} \snormo f_X \, \} \\
\noalign{\noindent and the upper Boyd index of $X$\ is defined by}
   q_X &=& \inf\{\, q : \exists c \ \forall a>1 \ \snormo{D_a f}_X \le c \,
   a^{-1/q} \snormo f_X \, \} .
\end{eqnarray*}
Thus we see that $1 \le p_X \le q_X \le \infty$.  Also, if $X$\ is the
Lorentz space $L_{p,q}$, then $p_X = q_X = p$.

We also define the {\it Hardy operators\/} as follows.
\begin{eqnarray*}
   H^{(p,r)} f(t) &=& 
   {1\over t^{1/p}} \left(\int_0^t (f^*(s))^r \, ds^{r/p} \right)^{1/r} ,\\
   H_{(q,r)} f(t) &=& 
   {1\over t^{1/q}}\left(\int_t^\infty (f^*(s))^r \, ds^{r/q}
   \right)^{1/r} ,\\
   H^{(p,\infty)} f(t) &=& 
   \sup_{0<s<t} (s/t)^{1/p} f^*(s) ,\\
   H_{(q,\infty)} f(t) &=& 
   \sup_{t<s} (s/t)^{1/q} f^*(s) ,\\
   H_{(\infty,r)} f(t)     &=& \left(\int_t^\infty {f^*(s)^r \over s} \, ds
   \right)^{1/r} ,\\
   H_{(\infty,\infty)} f(t) &=& f^*(t) .
\end{eqnarray*}

The Hardy operators play a very important role in interpolation theory.
The reason for this is the following result, essentially due to 
Holmstedt (1970).

\proclaim{\uppercase{Theorem 1}}  If $0 < p < q \le \infty$, and
$0<r,s \le \infty$, then
$$ \inf\{ t^{-1/p} \snormo{f'}_{p,r} + t^{-1/q} \snormo{f''}_{q,s} : 
   f' + f'' = f\}
   \approx
       H^{(p,r)} f(t) + H_{(q,s)} f(t)  . $$
\endproclaim
\proof In the case that $p=r$\ and $q=s$, this is the result in
Holmstedt (1970).  Otherwise, this result follows from Theorem~7 below.
\endproof

Boyd indices also play an important role in interpolation theory, because
the Boyd indices are strongly connected with the Hardy operators.  In fact,
the purpose of this paper is to make this connection firm.  In this
paper, we show the following result.  
The implications from left to right complement known results which would
yield the following in the case that 
$X$\ satisfied the triangle inequality (Maligranda, 1980, 1983).
The following result also generalizes a result from 
Ari\~no and Muckenhoupt (1990), where
they give necessary
and sufficient conditions for $H^{(1)}$\ to be bounded on a Lorentz
space.

\proclaim{\uppercase{Theorem 2}}  If $X$\ is a quasi-Banach r.i.\
space then we have the following.
\begin{itemlist}
\itemi for $0<p<\infty$\ and $0<r<\infty$\ the operator
$H^{(p,r)}$\ is bounded from $X$\ to $X$\ if and only if $p_X>p$.
\itemii For $0<q\le \infty$\ and $0<r<\infty$\ the operator
$H_{(q,r)}$\ is bounded from $X$\ to $X$\ if and only if $q_X<q$.
\itemiii for $0<p<\infty$\ the operator
$H^{(p,\infty)}$\ is bounded from $X$\ to $X$\ if $p_X>p$.
\itemiv For $0<q<\infty$\ the operator
$H_{(q,\infty)}$\ is bounded from $X$\ to $X$\ if $q_X<q$.
\end{itemlist}
\endproclaim

Note that the reverse implications are not true in parts~(iii) and~(iv).
For example, the operators $H^{(p,\infty)}$\ and $H_{(p,\infty)}$\ are
both bounded on the space $L_{p,\infty}$.

From this we can immediately generalize a result of Boyd (1967, 1969) 
to the following.

\proclaim{\uppercase{Theorem 3}}  
If $0 < p < q \le \infty$\ and $0 < r_1,r_2,s_1,s_2 
\le \infty$, 
and if $T:L_{p,r_1} \cap L_{q,s_1}
\to L_{p,r_2} \cap L_{q,s_2}$\ is a quasi-linear operator 
such that $\normo{Tf}_{p,r_1} \le c\, \normo f_{p,r_2}$\ and
$\normo{Tf}_{q,s_1} \le c\, \normo f_{q,s_2}$\ for all 
$f\in L_{p,r_1} \cap L_{q,s_1}$, and
if $X$\ is a quasi-Banach r.i.\ space with 
Boyd indices strictly between $p$\ and $q$, then
$\normo{Tf}_X \le c\, \normo f_X$\ for all $f\in L_{p,r_1} \cap L_{q,s_1}$.
\endproclaim
\proof
From Theorem~1, we see that
$$ H^{(p,r_1)}(Tf)(t) + H_{(q,s_1)}(Tf)(t) \le c 
   \bigl(H^{(p,r_2)}f(t) + H_{(q,s_2)}f(t)\bigr) .$$
Now the result
follows easily from Theorem~2 which implies that for $i=1,2$\ 
$$ \snormo{H^{(p,r_i)} + H_{(q,s_i)}}_X \approx \snormo f_X .$$
\endproof

Thus, as applications, we may obtain the following generalization of 
a result of Feh\'er, Gaspar and Johnen (1973).

\proclaim{\uppercase{Theorem 4}}  
The Hilbert transform is bounded on a quasi-Banach 
r.i.\ space $X$\ if and only if $p_X>1$\ and $q_X<\infty$.
\endproclaim
\proof
The implication from right to left follows immediately from Theorem~3.
As for the other way, this follows from the easy estimate:
$$ \hbox{\rm P.V.} \int_{-\infty}^\infty
   {f^*(y-x) \over y} \, dy \ge \textstyle{1\over 2}
   \bigl(H^{(1)}f(x) + H^*f(x)\bigr)  \qquad x > 0.$$
\endproof

We also obtain a result in the spirit of  
Ari\~no and Muckenhoupt (1990).  

\proclaim{\uppercase{Theorem 5}}  The Hardy--Littlewood maximal function is
bounded from $X(\R^n)$\ to $X(\R^n)$\ if and only if
$p_X>1$.
\endproclaim
\proof
Combine the argument given in Ari\~no and Muckenhoupt (1990) with Theorem~2 above.
\endproof

Before proceeding with the proof of Theorem~2, we will require the following
lemma.

\proclaim{\uppercase{Lemma 6}}  
Suppose that $X$\ is a quasi-Banach r.i.\ space.
Then given any $p>0$, there is a number $0<u\le p$\ such that
for any $f_1$, $f_2,\dots,$\ $f_n \in X$\ we have
$$ \normo{\left(\sum_{i=1}^n \modo{f_i}^p\right)^{1/p}}
   \le c \left(\sum_{i=1}^n \normo{f_i}^u\right)^{1/u} .$$
\endproclaim
\proof  Let $X^{(p)}$\ be the $p$-convexification of $X$, that is,
$ X^{(p)} = \{f: |f|^{1/p} \in X\}$\ and $\normo f_{X^{(p)}} = \normo{
|f|^{1/p}}^p$.  Clearly $X^{(p)}$\ is also a quasi-Banach space.  Thus
without loss of generality it is sufficient to show the above result
when $p = 1$.  But this follows immediately from 
Kalton, Peck and Roberts (1984), Lemma~1.1.
\endproof

\noindent
Proof of Theorem~2:
Let's consider the case for the lower Boyd indices.  The proof for the
other cases are almost identical.

Let us start by proving the implication from left to right.
Suppose that $p_X > p$.  Let us also suppose that $r < \infty$.  The
case when $r = \infty$\ then follows since $H^{(p,\infty)}f(t) \le
H^{(p,r)}f(t)$\ for any $0<r<\infty$.
Note that
\begin{eqnarray*}
   H^{(p,r)} f(t)
   &=&
   {1\over t^{1/p}}\left(\int_0^t (f^*(s))^r \, ds^{r/p} \right)^{1/r} \\
   &=&
   \left(\int_0^1 (D_af^*(t))^r \, da^{r/p}\right)^{1/r} \\
   &\le&
   \left(\sum_{n=-\infty}^0 2^{rn/p} (D_{2^n} f^*(t))^r \right)^{1/r} .
\end{eqnarray*}
Pick $0<u\le p$\ as given by Lemma~6.
Also, there is a number $p'\in (p,p_X)$\ such that
$$ \snormo{D_a f}_X \le c a^{-1/p'} \qquad (0<a\le 1). $$
Therefore,
\begin{eqnarray*}
   \snormo{H^{(p,r)} f}_X
   &\le&
   c \, \left(
   \sum_{n=-\infty}^0 2^{un/p} \snormo{D_{2^n} f}_X^u \right)^{1/u} \\
   &\le& c \, \left(
   \sum_{n=-\infty}^0 2^{un/p} 2^{-un/p'} \right)^{1/u} \snormo f_X \\
   &\le& c' \snormo f_X ,
\end{eqnarray*}
as desired.

Now let us prove the opposite implication.
Suppose that $\snormo{H^{(p,r)} f} \le C \, \snormo{f}$.
We are going to show that $X$\ has lower Boyd index greater than or equal
to $p/(1-1/C^r)$.  In order to do this, it is sufficient to show that
there is a number $0 < k < 1$\ such that for all numbers $a = k^n$\ for
integers $n \ge 1$, we have that 
$$ \snormo{D_a f} \le c\, a^{-(1-1/C^r)/p} \normo f .$$

Let us proceed.
By induction and a straightforward use of Fubini, we obtain the following
formula for the iteration:
$$ (H^{(p,r)})^{n+1} f(t) = {1\over t^{1/p}}\left( \int_0^t
   {(\log(\textstyle{t\over s})^{r/p})^n 
   \over n!} f^*(s)^r \, ds^{r/p} \right)^{1/r} ,$$
that is,
$$ (H^{(p)})^{n+1} f = \left(\int_0^1 
   {(\log\textstyle{1\over a^{r/p}})^n \over n!} 
   (D_a f^*)^p \, da^{r/p} \right)^{1/r}.$$
Note that ${(\log\textstyle{1\over a^{r/p}})^n 
\over n!} f^*(a)^r$\ is a decreasing
function in $a$, and hence for any $0<a<1$\ we have that
$$ (H^{(p)})^{n+1} f \ge \left(a^{r/p}
   {(\log\textstyle{1\over a^{r/p}})^n \over n!} \right)^{1/r} D_a f^* .$$
Hence
$$ \snormo{D_a f} \le C^{n+1} \left({n! \over a^{r/p}
   (\log\textstyle{1\over a^{r/p}})^n}
   \right)^{1/r} \snormo f .$$
Now let $k = \exp(-\textstyle{p\over r}C^r)$.  
Then using the estimate $n! \le c\, e^{-n} n^n$, 
we see that if $a = k^n$, then
$$ \snormo{D_a f} \le C c^{1/r} a^{-(1-1/C^r)/p} \normo f .$$
\endproof

Note that the proof actually gives a quite precise result.  For instance, it
is known that $\normo{H^{(1)} f}_p \le {p \over p-1} \normo f_p$.  The above
proof would show that the lower Boyd index is greater than or equal to
$p$, which is of course correct.

\section*{APPENDIX}

The result of this section is an extension of results already known in
interpolation theory.
Let $\normo f_{X+Y} = \inf\{\normo f'_X + \normo f''_Y:f'+f'' = f\}$,
and let 
$$ \normo f_{H(X,Y)} = \normo{f^*\big|_{[0,1]}}_X + 
   \normo{f^*\big|_{[1,\infty)}}_Y .$$

\proclaim{\uppercase{Theorem 7}}
Let $X$\ and $Y$\ be r.i.\ spaces such that
the following hold.
\begin{itemlist}
\itemi If $f$\ has support in $[0,1]$, then $\normo f_Y \le c_1 \normo f_X$.
\itemii If $f$\ is constant on intervals of the form $[n,n+1)$, then
$\normo f_X \le c_1 \normo f_Y$.
\itemiii $\normo{D_{1/4}f}_X \le c_2 \normo f_X$\ and
$\normo{D_{1/4}f}_Y \le c_2 \normo f_Y$.
\itemiv The quasi-triangle inequality constant for $X$\ and $Y$\ is less than
$c_3$.
\end{itemlist}
\noindent
Then
$$ \normo f_{X+Y} \le \normo f_{H(X,Y)} \le 2 c_1 c_2 c_3 \normo f_{X+Y} .$$
\endproclaim
\proof Clearly $\normo f_{X+Y} \le \normo f_{H(X,Y)}$.  To show the
opposite inequality, 
let
\begin{eqnarray*}
   E f(x) &=& \cases{
   f(x)            &  if $0\le x <1$ \cr
   0               &  if $x \ge 1 $,\cr} \\
   F f(x) &=& \cases{
   0               &  if $0\le x <1$ \cr
                   & \cr
   \displaystyle{ \int_{n}^{n+1} f(s) \, ds}
                   & if $n \le x < n+1$\ and $n$\ is
                     a positive integer.\cr} 
\end{eqnarray*}
Then we can see that $f^*(2x) \le (E + F) f^*(x) \le f^*(x/2)$.
Now suppose that
$ f^* = f_1 + f_2$.  Then 
$$ \normo{D_{1/2} f_1}_X
   \ge
   \normo{(E + F) f_1}_X
   \ge 
   {1\over 2} (\normo{E f_1}_X + \normo{F f_1}_X) 
   \ge
   {1\over 2c_1} (\normo{E f_1}_X + \normo{F f_1}_Y) ,$$
and
$$ \normo{D_{1/2} f_2}_Y
   \ge
   \normo{(E + F) f_2}_Y
   \ge 
   {1\over 2} (\normo{E f_2}_Y + \normo{F f_2}_Y) 
   \ge
   {1\over 2c_1} (\normo{E f_2}_X + \normo{F f_2}_Y) .$$
Hence
$$ \normo{D_{1/2} f_1}_X + \normo{D_{1/2} f_2}_Y
   \ge
   {1\over 2c_1 c_3} (\normo{E f^*}_X + \normo{F f^*}_Y) ,$$
and so $\normo{(E+F)f^*}_{H(X,Y)} \le 2c_1 c_3
\normo{D_{1/2} f^*}_{X+Y}$.  Therefore,
$$ \normo f_{H(X,Y)} 
   \le \normo{(E + F)D_{1/2} f^*}_{H(X,Y)}
   \le 2c_1 c_3 \normo{D_{1/4} f}_{X+Y} 
   \le 2 c_1 c_2 c_3 \normo f_{X+Y} . $$
\endproof

\end{document}